\magnification1100
\hsize15.6truecm
\vsize22.8truecm
\hoffset-.13truecm

\def\Sym{ {\cal S}ym }
\def\puce{\vskip 5mm\noindent {$\scriptstyle\bullet$}\kern 8pt}
\def\d{\partial}

\overfullrule=0pt
\baselineskip 14pt
\parindent=.5truein
\hfuzz=3.44182pt

\font\title=cmbx10 scaled\magstep2
\font\normal=cmr10 
\font\ninerm=cmr9 

\font\germten=eufm10
\def\Sn{\hbox{\germten S}}
\font\FontSetTen=msbm10 
\def\N{\hbox{\FontSetTen N}}

\def\R{\hbox{\FontSetTen R}}

\def\Q{\hbox{\FontSetTen Q}}

\font\titlefont=cmbx10 scaled\magstep1
\font\bol=cmbx8 scaled\magstep2
 
\def\la {\lambda}
\def\d {\partial}

\def \sas {\vskip .06truein}
\def\sa{{\vskip .125truein}}

\def\sap{{\vskip .25truein}}

\def \ses {\enskip = \enskip}
\def \s=s {\enskip = \enskip}
\def \sps {\enskip + \enskip}

\def \ssp {\hskip .25em}

\def \sqr#1#2{{\vcenter{\vbox{\hrule height.#2pt
	\hbox{\vrule width.#2pt height#1pt \kern#1pt
		\vrule width.#2pt}
		\hrule height.#2pt}}}}
\def \square{\mathchoice\sqr68\sqr68\sqr{2.1}3\sqr{1.5}3}
\vsize=8truein
\sap
\def\today{\ifcase\month\or
January\or February\or March\or April\or may\or June\or
July\or August\or September\or October\or November\or
December\fi
\space\number\day, \number\year}

\normal
\centerline{ \titlefont Determinantal expressions for Macdonald polynomials}
\sap 
\centerline{\bf L. Lapointe, A. Lascoux, and J. Morse}
\sap
\centerline{\it Institut Gaspard Monge, C.N.R.S.} 
\centerline{\it Universit\'e de Marne-la-Vall\'ee, 5 Bd Descartes,}
\centerline{\it Champs sur Marne, 77454 Marne La Vall\'ee Cedex, France}
\sap
{
\ninerm
\smallskip\midinsert\narrower\noindent\baselineskip11.5pt
\noindent{ABSTRACT:}
We show that the action of classical operators associated to the
Macdonald polynomials on the basis of Schur functions, 
$S_{\lambda}\left[ X (t-1)/(q-1) \right]$, 
can be reduced to addition in $\lambda-$rings.
This provides explicit formulas for the Macdonald polynomials
expanded in this basis as well as in the ordinary Schur basis,
$S_{\lambda}[X]$, and the monomial basis, $m_{\lambda}[X]$.
\endinsert
}
\sap
\noindent {\bol Introduction }
\sa
Important developments in the theory of symmetric functions
rely on the use of the Macdonald polynomial basis, $\{J_\la(X;q,t)\}_\la$ [M1].
This basis specializes to several fundamental bases including the
Schur, Hall-Littlewood, Zonal, and Jack.
It has been conjectured that the Macdonald polynomials 
occur naturally in representation theory of the 
symmetric group [GH].  It is also known that these polynomials are 
eigenfunctions of a family of commuting
difference operators with significance in many-body physics [RS]. 
\sas
We first reformulate the Macdonald operators in terms 
of divided differences, operators that act naturally 
on the Schur function basis.
This enables us to show that the action of the Macdonald 
operators on the modified Schur basis, $S_\lambda[X(t-1)/(q-1)]$, 
primarily amounts to addition in $\lambda$-rings (Theorem 2.1).
The action provides a determinantal expression for the 
Macdonald polynomials expanded in this Schur basis.
By involution this expression can be converted
into an explicit formula for the Macdonald polynomials
in the usual Schur basis, $S_\lambda[X]$, and in the monomial basis,
$m_{\lambda}[X]$ (Theorems 3.1 and 3.3).
Further, it is known that the Macdonald polynomials can be built
recursively using `creation operators'.
We reformulate these operators in a manner similar to our 
expression for the Macdonald operators enabling us to give a simple
expression for their action on the modified Schur functions (Theorem 4.2).
\sap
\noindent {\bol 1. Background}
\sa
To begin, we clarify some notation. 
$\Sym$ refers to the ring of symmetric functions.
A partition will be defined as a weakly decreasing sequence
of non-negative integers,
$\lambda=(\lambda_1,\ldots,\lambda_n)$
with $\ell(\lambda)$, the number of parts of $\lambda$ and
$\lambda'$ the conjugate to $\lambda$. 
To a partition $\lambda$, there
corresponds a Ferrers' diagram with
$\lambda_i$ lattice squares in the $i^{th}$
row, from the bottom to top.  For each square $s$ in the diagram of $\lambda$,
we define $a(s)$ and $l(s)$ to be the number of squares respectively to the
north and east of $s$.  The {\it order of $\lambda$} is the sum of the parts of
$\lambda$, denoted $|\lambda|$.  
The {\it dominance order} on $\lambda$ is defined such
that for $|\lambda|=|\mu|$, 
$\lambda\leq\mu$ when $\lambda_1+\cdots+\lambda_i\leq
\mu_1+\cdots+\mu_i$ for all $i$.  
\sas
We shall use $\lambda$-rings, needing only the formal ring of symmetric
functions $\Sym$ to act on the ring of rational functions in $x_1,\dots,x_n,q,t$,
with coefficients in $\R$.  The ring $\Sym$ is generated by power sums $\Psi_i$, $i=1,2,3\dots$.
The action of $\Psi_i$ on a rational function $\sum_{\alpha} c_{\alpha} u_{\alpha}/ \sum_{\beta} d_{\beta} v_{\beta}$
is by definition
$$
\Psi_{i} \left[ {\sum_{\alpha} c_{\alpha} u_{\alpha}  \over \sum_{\beta} d_{\beta} v_{\beta} } \right]
 ={\sum_{\alpha} c_{\alpha} u_{\alpha}^i \over \sum_{\beta} d_{\beta} v_{\beta}^i}, \eqno 1.1
$$
with $c_{\alpha},d_{\beta} \in \R$ and $u_{\alpha}, v_{\beta}$ monomials in $x_1,\dots,x_n,q,t$.  Since any symmetric
function is uniquely expressed in terms of the power sums, formula 1.1 extends to an action of
$\Sym$ on rational functions.  In particular, a symmetric function 
$f(X)$ of $X=\{x_1,\dots,x_n \}$ can be denoted $f[x_1+\cdots+x_n]$.
\sas
The Schur function $S_{\lambda}$, with $Y, Z \in \Q [x_1,\dots,x_n,q,t]$, 
is such that
$$
S_{\lambda}[Y \pm Z] \ses \sum_{\mu} S_{\lambda/\mu}[Y] S_{\mu}[\pm Z]\ssp, \eqno 1.2
$$
where $S_{\mu}[-Z]= (-1)^{|\mu|} S_{\mu'}[Z]$, with $\mu'$ the 
partition conjugate to $\mu$.  
Schur functions will be considered as determinants of complete functions:
$$
S_{\lambda}[X] = {\hbox {det} }\left|
\matrix{ S_{\lambda_1} [X]  & S_{\lambda_1+1}[X]  & \cdots & S_{\lambda_1+n-1}[X]  \cr
  S_{\lambda_2 -1}[X]  & S_{\lambda_2}[X]  & \cdots & S_{\lambda_2 +n-2}[X]  \cr
  \vdots & \vdots & \ddots & \vdots \cr
  S_{\lambda_n -n+1}[X]  & S_{\lambda_n -n+2}[X]  & \cdots & S_{\lambda_n}[X] \cr} 
  \right| , \eqno 1.3
$$
as will skew Schur functions:
$$
S_{\lambda/\mu}[X] = {\hbox {det} }\left|
\matrix{ S_{\lambda_1-\mu_1} [X]  & S_{\lambda_1+1 -\mu_2}[X]  & \cdots & S_{\lambda_1+n-1 -\mu_n}[X]  \cr
  S_{\lambda_2 -1 -\mu_1}[X]  & S_{\lambda_2- \mu_2}[X]  & \cdots & S_{\lambda_2 +n-2 -\mu_n}[X]  \cr
  \vdots & \vdots & \ddots & \vdots \cr
  S_{\lambda_n -n+1 -\mu_1}[X]  & S_{\lambda_n -n+2 -\mu_2}[X]  & \cdots & S_{\lambda_n -\mu_n}[X] \cr} 
  \right| \ssp.
\eqno 1.4
$$
In these expressions, $S_{i}=0$ for $i<0$.
\sas
We denote $X=x_1+\cdots+x_n$ and a partial alphabet, $X_I=x_{i_1}+\cdots+x_{i_k}$, for $I = \{i_1, \dots,i_k \}$.
The complement of $X_I$ in $X$ will be $X_I^c$.   
For $I=\{1, \dots,k \}$, we have $X_k= x_1 + \cdots +x_k$.
Two particular elements, $X^{tq}:=X(t-1)/(q-1)$ and $X^t:=X(t-1)$ 
will be used throughout the paper.
\sas
We shall need the divided differences $\partial_{\mu}$, indexed by elements $\mu$ of the
symmetric group $\Sn (n)$.  
In particular, if $\mu$ is the simple transposition $\sigma_i$, we have
$$\partial_i  \, : f   \longrightarrow (f-\sigma_i f)/ (x_i-x_{i+1}) \, . \eqno 1.5
$$
The divided differences satisfy the 
Coxeter relations [M2]
$$\d_i \d_j  = \d_j \d_i, \quad |i-j| > 1 ; \qquad \quad 
  \d_i \d_{i+1} \d_i  = \d_{i+1} \d_i \d_{i+1} \ssp,             \eqno 1.6
  $$
which imply that $\partial_{\mu}$ can be written 
as the product of elementary operators $\partial_i$ 
corresponding to any reduced decomposition of $\mu$.
It should be noted that symmetric functions in $x_i,x_{i+1}$ are scalar with
respect to $\partial_i$  and that $\d_{i}  \cdot 1 =0$.  Consequently,
symmetric functions in $x_1,\ldots,x_n$ commute with all $\partial_\mu$.
\sas
We have that $\partial_{n-1}\cdots\partial_{1}$ is the {\it Lagrange operator}
[L1],
$$
\partial_{n-1}\cdots\partial_{1}: f\in \Sym(1|n-1) \longrightarrow \sum_i 
{f(x_i,x_i^c)\over 
R(x_i,x_i^c)}
\in \Sym (n)\ssp, \eqno 1.7
$$
or more generally, 
$\partial_{(k|n-k)} =
(\partial_{n-k}\cdots\partial_{1})\cdots
(\partial_{n-2}\cdots\partial_{k-1})
(\partial_{n-1}\cdots\partial_{k})$ is the {\it Sylvester operator} [L1], 
$$
\partial_{(k|n-k)}:  f\in \Sym(k|n-k) \longrightarrow \sum_{|I|=k} 
{f(X_I, X_I^c)\over R(X_I,X_I^c)}
\in \Sym(n)\ssp,
\eqno 1.8
$$
where $\Sym(k|n-k)$ is the space of functions symmetrical in $x_1,\dots,x_k$ and in $x_{k+1},\dots,x_n$, 
and 
where $R(X,Y)$ is the resultant of two alphabets $X$ and $Y$ ;
$$
R(X,Y)\ses
\prod\nolimits_{x\in X, y\in Y}
(x-y)\ssp. 
\eqno 1.9
$$
\sas
These operators may be used to define the Euler-Poincar\'e 
characteristic $\chi$.  In fact, in his study
of Riemann-Roch theorem, Hirzebruch defined a `$\chi_y$'-characteristic, 
or `$\chi_y$'-genus [H],[HBJ].  We shall use the one corresponding 
to a relative flag manifold, changing $y$ into $1/t$ .
We define the following operators on polynomials in $x_1,\dots,x_n$,
$$
\chi_{(1|n-k)}^{(k)}\ses
\partial_{n-1}\cdots\partial_{k}\,
R(x_k,X_k^c/t) \, ,
\eqno 1.10
$$
considered as a composition of multiplication 
by the resultant followed by a sequence of divided differences.
Note that the image of 1 is
$(1+\cdots +t^{-n+k})$ and that the superscript refers
to the variable at which the symmetrization begins.
We shall also use 
$$
\chi_{(k|n-k)}
\ses
\partial_{(k|n-k)}
\, R(X_k, X_k^c/t)
\eqno 1.11
$$
and
$$
\chi_{\omega(k)}\ses
\partial_{\omega(k)}\,
\prod_{1 \leq i<j \leq k } (x_i- x_j/t)
\ssp,
\eqno 1.12
$$
where $\omega(k)$ is the longest permutation in $\Sn(k)$
and $\chi_{\omega(k)}\cdot 1= k_t!= t^{-k(k-1)/2} (t)_k/(1-t)^k\,,$ where
$(t)_k=(1-t)\cdots(1-t^k)$.
To compute with $\d_{\omega(k)}$, we shall use the factorization
$\partial_{\omega(k)}=
(\partial_{k-1}) (\partial_{k-2}\partial_{k-1}) 
\cdots (\partial_{1}\cdots\partial_{k-1}).$
\sas
When operating on elements of $\Sym(\ell|k-\ell)$, one has
$$
\chi_{\omega(k)}\,f\ses
\ell_t!\,(k-\ell)_t!\,
\chi_{(\ell|k-\ell)}\,f
\quad\quad\forall\,\, f\in Sym(\ell|k-\ell)\ssp.
\eqno 1.13
$$
We shall also need the factorization
$$
\chi_{(k|n-k)}
\chi_{(1|k-1)}^{(1)}\cdots
\chi_{(1|2)}^{(k-2)}
\chi_{(1|1)}^{(k-1)}\,
\ses
\chi_{(1|n-1)}^{(1)}\cdots
\chi_{(1|n-k+1)}^{(k-1)}
\chi_{(1|n-k)}^{(k)}\, 
\eqno 1.14
$$
and the following lemma (proved in the appendix).
\sa
\noindent{\bf Lemma 1.1}\quad {\it For an alphabet $C$, symmetric in $x_1,\dots,x_n$,
and $c(k)= t^{-1}+\cdots+t^{-n+k}$ , one has}
$$
\chi_{(1|n-k)}^{(k)}\,
S_j[C+X_k^t]
\ses
c(k)\,S_j[C+X_{k-1}^t]+S_j[C+X^t] \, .
\eqno 1.15
$$
\sap
\noindent {\bol 2. Macdonald Operators}
\sa
This section presents the action of
the Macdonald operators on $S_{\lambda}[X^{tq}]$ 
as a sum of determinants.
The Macdonald operators [M1] are defined as
$$ M_k\ses 
t^{k \choose 2}  \sum_{|I|=k}
\,R(tX_I,X_I^c)\,
\Theta_{i_k} \cdots \Theta_{i_1} 
\qquad \qquad k=1,2,3,\dots \, ,  \eqno 2.1
$$
where the sum is over all $k$-element subsets $I=(i_1,\dots,i_k)$ of 
$\{ 1, \dots, n\}$ and 
$$\Theta_l:f[X]\longrightarrow f[X+(q-1)x_l]
\ssp.
\eqno 2.2
$$
It is more convenient to use a reformulation for these operators in terms
of divided differences,
$$ 
\eqalign{M_k & \ses t^{{k \choose 2}+ k(n-k)} \chi_{(k|n-k)} \Theta_k \cdots \Theta_1 \cr
& \ses 
t^{k \choose 2}\,
\d_{(k|n-k)} \,
R (t X_k, X_k^c) \,
\Theta_k \cdots \Theta_1 
\ssp, \cr}
\eqno 2.3
$$
thanks to 1.8.

We shall express the action of $M_k$ on $S_\lambda[X^{tq}]$ with formal
operators,
$\left[ \matrix{ a_1 & a_2 & \cdots& a_n \cr
   b_1 &b_2 & \cdots &b_n }\right]$, acting on the
columns of an $n \times n$ matrix of Schur functions
such that for all elements in any column $\ell$ , 
$S_{\beta}[Y]\rightarrow S_{\beta+a_{\ell}}[Y+b_{\ell}]$.  
For example;
$$
\left[ \matrix{ 0 & 30 & -5 \cr
   b_1 & Y & 0 }\right]\,:\,
\left( \matrix{ 
S_{11}[A] & S_{4}[B]& S_{8}[C] \cr
S_{0}[A] & S_{2}[B]& S_{13}[C] \cr
S_{1}[A] & S_{2}[B]& S_{3}[C] \cr} 
 \right)  \rightarrow  
\left( \matrix{ 
S_{11}[A + b_{1}] & S_{34}[B+Y]& S_{3}[C] \cr 
S_{0}[A+ b_{1}] & S_{32}[B+Y]& S_{8}[C] \cr 
S_{1}[A+ b_{1}] & S_{32}[B+Y]& S_{-2}[C] }
\right)\ssp . \eqno 2.4
$$
\sa
\noindent {\bf Theorem 2.1}\quad  {\it The action of the Macdonald
operators on $S_{\lambda}[X^{tq}]$, with $\ell(\lambda) \leq n$, can be expressed as}
$$
M_k S_{\lambda}[X^{tq}]\ses 
\sum_{1 \leq i_1 < \cdots
 < i_k \leq n} 
t^{(n-i_1)+ \cdots+ (n-i_k) } \, 
\left[ \matrix{ 0 & \cdots & 0 &\cdots&0 &\cdots &  0\cr
   0 &\cdots &X^t &\cdots&X^t& \cdots &0 }\right] \, S_{\lambda}[X^{tq}]
\ssp,
\eqno 2.5
$$
{\it where $X^{t}$ occurs in columns $i_1,\dots,i_k$.}  
\sa
\noindent For example, if $n=2$, 
$$
M_1 S_{(5,2)} \ses
t \left|\matrix{S_5[X^{tq}+X^t]\hfill &S_6 [X^{tq}]\hfill\cr
S_1[X^{tq}+X^t]\hfill &S_2[X^{tq}]\hfill}
\right|\sps
\left|\matrix{S_5[X^{tq}] \hfill &S_6[X^{tq}+X^t]\hfill\cr
S_1[X^{tq}]\hfill &S_2[X^{tq}+X^t]\hfill}
\right| \ssp, \eqno 2.6
$$
or more concisely since $X^{tq}+X^t =qX^{tq}$, using $S_i$ for $S_i[X^{tq}]$ ,
$$
M_1 S_{(5,2)} \ses
t \left|\matrix{q^5S_5\hfill &S_6\hfill\cr
qS_1\hfill &S_2\hfill}
\right|\sps
\left|\matrix{S_5\hfill &q^6S_6\hfill\cr
S_1\hfill &q^2S_2\hfill}
\right| \ssp. \eqno 2.7
$$
\sa
\noindent {\bf Proof of Theorem 2.1} \quad  
The Macdonald operators, 
$$ 
M_k \ses 
t^{k \choose 2}\,
\d_{(k|n-k)} \,
R (t X_k, X_k^c) \,
\Theta_k \cdots \Theta_1 \ssp, \eqno 2.8
$$
act first with a product of $\Theta_i$,
defined to replace $x_i$ by $qx_i$. 
However, this $q$-deformation is actually a $t$-deformation when acting on $f[X^{tq}]$ since
$$
\Theta_{\ell}:f[X^{tq}]\rightarrow f \left[ \bigl( X + (q-1) x_{\ell} \bigr) \cdot {t-1\over q-1} \right]=f[X^{tq}+(t-1)x_{\ell}] 
\ssp. \eqno 2.9
$$
Now for $Y=X^{tq}-X$, if we define 
a formal operator $\Theta_{\ell}^t$ 
that sends $x_i\rightarrow x_i$ and
$x_{\ell}\rightarrow tx_{\ell}$ while leaving $Y$ invariant;
i.e. $Y+X\rightarrow Y+X+(t-1)x_{\ell}$, we obtain  
$$
R(tX_k,X_k^c)\,\Theta_k\cdots\Theta_1\,S_\lambda[Y+X]
\ses
\Theta_k^t\cdots\Theta_1^t
\,R(X_k,X_k^c)\,S_\lambda[Y+X]\ssp, \eqno 2.10
$$
where $\Theta_i$ acts on both $Y$ and $X$ on the 
left hand side of this expression.
Further, we have the following determinantal expression
for $R(X_k,X_k^c)\,S_\lambda[Y+X]$;
$$ 
\eqalign{
R(X_k,X_k^c)&\, S_\lambda[Y+X] \ses\cr
& \Big|S_{\alpha_i-k+n}[Y+X_k] \cdots 
S_{\alpha_i-1+n}[Y+X_k] 
S_{\alpha_i+k}[Y+X_k^c]
\cdots S_{\alpha_i+n-1}[Y+X_k^c]
\Big|_{i=1\cdots n}
\ssp,
}
\eqno 2.11
$$
where $\alpha_i\rightarrow \lambda_{i}-i+1$
(to obtain this expression see the appendix). 
Since $\Theta_k^t\cdots\Theta_1^t$ acts exclusively 
on the first $k$ columns of this determinant, we have that 
$$ \eqalign{ &t^{-{k \choose 2}} M_kS_{\lambda}[Y+X] \ses \cr
&   \d_{(k|n-k)}
  \Big|S_{\alpha_i-k+n}[Y+tX_k] \cdots
S_{\alpha_i-1+n}[Y+tX_k] 
S_{\alpha_i+k}[Y+X_k^c]
\cdots S_{\alpha_i+n-1}[Y+X_k^c]\Big|_{i=1..n}
\ssp. \cr}
\eqno 2.12
$$
As it remains to show that this is equivalent
to expression 2.5, we proceed by giving the action of 
$\d_{(k|n-k)}$ on products of arbitrary 
minors of order $k$ on the first $k$ columns and of order 
$n-k$ on the last $n-k$ columns of this determinant. 
The action of $\d_{(k|n-k)}$ on the whole determinant will then
follow by linearity.
\sa
{\it{
\noindent{\bf Lemma 2.2} \quad For $\alpha=(\alpha_1,\dots,\alpha_k)$,
$\beta=(\beta_1,\dots,\beta_{n-k})$, and $A, B$ 
invariant under $\Sn (n)$, we have
$$
\eqalign{ \d_{(k|n-k)}  \Bigl( & S_{\alpha}[A+tX_k]  S_{\beta}[B+X_k^c] \Bigr) 
  \ses t^{k(n-k) }\sum_{I} (-1/t)^{|I|} \cr
& \quad \times \left(  \left[
\matrix{-n+k+i_k & \cdots & -n+k+i_1 \cr
 tX & \cdots & tX \cr} \right]S_{\alpha}[A]\right) \,\left( \left[ \matrix{-i_1' & \cdots & -i_{n-k}' \cr
 X & \cdots & X \cr} \right] S_{\beta}[B] \right)
\, ,\cr} \eqno 2.13
$$
where $I=(i_1,\dots,i_k)$ and $(i_1',\dots,i_{n-k}')$ are pairs
of conjugate partitions.}}
\sa
\noindent {\bf Proof of lemma} \quad 
Define the elements $A'=A+tX$ and $B'=B+X$.  
The product  $S_{\alpha}[A+tX_k]  S_{\beta}[B+X_k^c]$ can be
expanded, thanks to 1.2, into products of Schur functions 
in $-X_k^c$ and $X_k$ (with coefficients in $A'$ and $B'$), 
$$
\eqalign{
S_{\alpha}[A+tX_k] \, S_{\beta}[B+X_k^c] & \ses S_{\alpha}[A'-tX_k^c] \, S_{\beta}[B'-X_k] \cr
 & \ses \sum_{I,J} S_{\alpha/I}[A'] t^{|I|} S_I [-X_k^c] S_{\beta/J} [B'] (-1)^{|J|} S_{J'}[X_k] \, ,
\cr}
\eqno 2.14
$$
where $I$ and $J$ are partitions and $J'$ is conjugate to $J$.  
Since $\ell(I),\ell(J') \leq k$ and $\ell(I'),\ell(J) \leq n-k$, where
$I'$ is conjugate to $I$,  we have  $I,J' \subseteq (n-k,\dots,n-k)$.
Thus the following property describing two adjoint bases of $\Sym (k,n-k)$ 
as a free module over $\Sym (n)$ ([L2], Corollary 2.2),
$$
\d_{(k|n-k)} \Bigl( S_I [-X_k^c] \, S_J [X_k] \Bigr) = \cases{ 1 &  if  $ I=(n-k-J_k,\dots,n-k-J_{1})$ \cr
         0 &  else}  \ssp,
\eqno 2.15
$$
gives
$$
\d_{(k|n-k)}\Bigl(S_{\alpha}[A+tX_k] \, S_{\beta}[B+X_k^c] \Bigr) = 
t^{k(n-k) }\sum_{J} (-1/t)^{|J|} S_{\alpha/(n-k-J_k',\dots,n-k-J_{1}')}[A']  S_{\beta/J} [B'] .
\eqno 2.16
$$
The lemma is proved by letting $J\to I'$ in the last expression and using the definition 1.4 of skew Schur functions. \hfill  $\square$
\sa
Now applying Lemma 2.2 directly to expression 2.12, we have
$$
\eqalign{
&  M_kS_\lambda [Y+X] \cr
& \ses
\sum_I (-1/t)^{|I|} t^{k(n-k)+k(k-1)/2}  \,
\left[ \matrix{-n+k+i_k & \cdots &-n+k+i_1 & -i_{1}' &\cdots & -i_{n-k}' \cr
   tX &\cdots &tX &X &\cdots &X }\right] \cr
& \qquad \qquad \qquad \qquad \qquad \times  \Big| S_{\alpha_i-k+n}[Y] \cdots S_{\alpha_i-1+n}[Y] S_{\alpha_i+k}[Y] 
\cdots S_{\alpha_i+n-1}[Y] \Big|_{i=1..n} \cr\cr
& \ses 
\sum_I (-1/t)^{|I|} t^{k(n-k)+k(k-1)/2} \,
\left[ \matrix{ i_k & \cdots & i_1+k-1 & -i_{1}'+k &\cdots & -i_{n-k}'+n-1\cr
   tX &\cdots &tX &X &\cdots &X }\right]\cr
&\qquad\qquad\qquad \qquad \qquad \times
\Big| S_{\alpha_i} [Y] \cdots S_{\alpha_i} [Y] S_{\alpha_i} [Y] 
\cdots S_{\alpha_i}[Y] \Big|_{i=1..n}
\ssp.
}
\eqno 2.17
$$
\noindent For $\rho=(0,\ldots,n-1)$ and $\lambda^\omega$ the reverse reading of
a partition $\lambda$, we have from [M1] that
$\rho + \lambda^\omega \cup -\lambda'= \sigma\rho$ 
for some permutation $\sigma$. The top row
of the formal operator in 2.17, $\rho+(i_k,\ldots,i_1,-i_1',\ldots,-i_{n-k}')$, may 
thus be rearranged to $\rho$.  The rearrangement sends
column $\ell+1$, containing $tX$ and corresponding to index 
$i_{k-\ell}+\ell$, to column $i_{k-\ell}+\ell+1$.
This requires $i_{k-\ell}$ commutations with $X$-columns 
inducing a $(-1)^{i_{k-\ell}}$.  
Proceeding on all $tX$-columns yields
$$ 
M_k \,S_\lambda [Y+X]
 \ses \sum_I t^{k(n-k) +k(k-1)/2-|I|} \,
\left[ \matrix{ 0 & \cdots & 0 &\cdots & 0 &\cdots &  0\cr
   X &\cdots &tX & \cdots& tX &\cdots &X }\right]
S_\lambda [Y]  \ssp, \eqno 2.18
$$
where the $tX$ are in columns $i_k+1,\dots,i_1+k$.  
Finally, the substitution $Y= X^{tq}-X$ with the shift $i_k+1\rightarrow i_k,\ldots,
i_1+k\rightarrow i_1$ gives Theorem 2.1 .\hfill $\square$
\sa
\noindent Remark:  The image of $S_{\lambda}[X^{tq}]$ under
the Macdonald operator $M_k$ in the basis of products of complete functions
$S^{\mu}[X^{tq}]=S_{\mu_1}[X^{tq}] S_{\mu_2}[X^{tq}] \dots$ can
be derived from 2.5.  However, there exists a scalar product for which
the $S^{\mu}[X^{tq}]$ basis is adjoint to the $m_{\mu}[X]$ basis,
the $S_{\lambda}[X^{tq}]$ basis is adjoint to the $S_{\lambda'}[X]$ basis, and 
$M_k$ is self-adjoint.  This allows that formula 2.5 be connected 
to the action of the Macdonald operators on monomial functions given in 
[M1],VI 3.6.
\sap
\vfill\supereject
\noindent {\bol 3. Macdonald Polynomials}
\sa
This section contains determinantal expressions for 
the Macdonald polynomials in the bases,
$S_{\mu}[X^{tq}]\, , S_{\mu}[X]\,$, and 
$m_\mu[X]\,$, obtained from the action of the 
operator $M_1$ on $S_\mu[X^{tq}]$.  
\sas
The Macdonald polynomials [M1] can be defined to be the eigenfunctions of 
the Macdonald operators.  In particular we have , for any partition $\lambda$, 
$$ 
M_1\,J_{\lambda}(X;q,t) \ses [|\lambda|]  \,J_{\lambda}(X;q,t)\ssp,
\eqno 3.1
$$
where the symbol
$$
[|\alpha|] \, := \, q^{\alpha_1}t^{n-1}+q^{\alpha_2}t^{n-2} + \cdots+ q^{\alpha_n}
\ssp, \eqno 3.2
$$
will be used for any composition $\alpha\in\N^n$.
Since the eigenvalues are distinct for each $\lambda$, 
any polynomial that satisfies 3.1 must be
proportional to $J_{\lambda}(X;q,t)$.  It is with this in mind that
we proceed to construct such polynomials using the action of 
$M_1$ on $S_\lambda[X^{tq}]$.
\sas
First, we have that the ordered expansion of $S_{\mu}$ in terms
of complete functions is given, for $S^\alpha=S_{\alpha_1}\cdots S_{\alpha_n}$, by
$$
S_{\mu}[X] 
\ses 
\sum_{\sigma \in \Sn (n)} (-1)^{\ell(\sigma)} 
\, 
S^{\sigma(\mu+\rho)  -\rho}[X] 
\, := \,
\sum_{\alpha } \epsilon(\mu,\alpha)
\, 
S^{\alpha}[X] 
\ssp, \eqno 3.3
$$
summed over all vectors $\alpha \in  \N^{n}$, 
where $\epsilon(\mu,\alpha) \in \{ 0,+1,-1\}$.
The set of all $\alpha$ with a corresponding nonzero
$\epsilon(\mu,\alpha)$ is included in 
the set
$\{\bar \sigma \mu \, |\, \bar \sigma \in \Sn ({\ell(\mu)}) \}$
where we define
$ \bar \sigma_i \mu= -(\dots,\mu_{i+1}-1,\mu_i+1,\dots)$.  
For example, if $n=3$,
$$
\matrix{ &&  &-(0,4,1)&{\bar\sigma_2\atop\longrightarrow} &(0,0,5)& \cr
	   &&{\bar\sigma_1\atop \nearrow}& & &&{\bar\sigma_1\atop\searrow}\cr
\mu=(3,1,1)\,:\qquad& (3,1,1)&& & & & &-(-1,1,5)\cr
	   &&{\searrow\atop\bar\sigma_2} & & &&{\nearrow\atop\bar\sigma_2} \cr
           &&  &-(3,0,2)&{\bar\sigma_1\atop\longrightarrow} &(-1,4,2)& } 
\ssp.
\eqno 3.4
$$
For $\mu = (3,1,1)$, 
$\epsilon(\mu,\alpha)\,\alpha = 
\{(3,1,1),-(0,4,1),(0,0,5),-(3,0,2)\}$,
since $\alpha\in \N^n$ cannot have negative components. 
\sas
Now, recall that Theorem 2.1 gives  explicitly
$$
M_1\,S_{\mu} [X^{tq}]
\ses 
\sum_{1 \leq i \leq n} 
t^{n-i} \, \left[ \matrix{ 0 & \cdots &0 &\cdots &  0\cr
   0  &\cdots&X^t& \cdots &0 }\right]
S_{\mu} [X^{tq}]\ssp. \eqno 3.5
$$
The expansion of the right hand side,
taking into account that $S_k[X^{tq}+X^t]=q^kS_k[X^{tq}]$,
is
$$
M_1\,S_{\mu}[X^{tq}] 
=
\sum_{\sigma \in \Sn (n)} (-1)^{\ell(\sigma)} 
\, [| \sigma (\mu+\rho)  -\rho |] \,
S^{\sigma (\mu+\rho)  -\rho}[X^{tq}] 
=
\sum_{\alpha } \epsilon(\mu,\alpha)
\, [|\alpha|]  \,
S^{\alpha}[X^{tq}] 
\ssp.\eqno 3.6
$$
The action of the Macdonald operator $M_1$ is triangular on the basis,
$S_\lambda[X^{tq}]$, giving that
$$
M_1\,S_{\mu}[X^{tq}] 
\ses 
[|\mu|] \,S^{\mu}[X^{tq}] \,+
\sum_{\alpha} \epsilon(\mu,\alpha) \,
S^{\nu}[X^{tq}] 
\ssp,
\eqno 3.7
$$
where the sum is over all $\alpha$ a permutation of partitions $\nu > \mu$.  
Formula 3.7 demonstrates that the eigenspaces of $M_1$ are $1$-dimensional 
and reveals the eigenvalues.
\sas
For any pair of partitions $\mu,\nu \in \N^n$ of the same weight, let us write
$$
[\lambda]_{\mu\nu} (q,t) \ses
 \sum_{\alpha } \epsilon(\mu,\alpha) \,
\bigl( [|\lambda|] -[|\alpha|]  \bigr)\ssp, \eqno 3.8 $$
where $\alpha$ runs over all distinct permutations of $\nu \in \N^n.$
Note that, in the following,  $[\lambda]_{\mu \nu}$ shall stand for $[\lambda]_{\mu \nu}(q,t)$.
\sa
\noindent {\bf Theorem 3.1}
\quad 
{\it Expanded in terms of $S_\lambda[X^{tq}]$,
the Macdonald polynomials are 
$$ J_{\lambda}(X;q,t) = {c_{\lambda'}(t,q) \over v_{\lambda}(q,t)} \det \left| \, 
\matrix{ S_{\lambda} [X^{tq}] & \ldots & S_{\mu}[X^{tq}] & \ldots \cr
                                                \vdots& & \vdots &  \cr
 [\lambda]_{\lambda\nu}  & \ldots &  
 [\lambda]_{\mu\nu}  & \ldots \cr
 \vdots & & \vdots &  \cr} \right|_{\nu>\lambda\atop\mu\geq\lambda}\ssp. 
\eqno 3.9
$$
Columns are indexed by all partitions $\mu \geq \lambda$.  Aside from the first row,
the entries are the polynomials in $q,t$, $[\lambda]_{\mu \nu}$, $\nu > \lambda$.
The normalization factor consists of 
$$ v_{\lambda}(q,t) = \prod_{\mu > \lambda} 
\Bigl( [|\lambda|]  -[|\mu|]  \Bigr)= \prod_{\mu > \lambda} 
 [\lambda]_{\mu \mu}  \ssp, \eqno 3.10$$
and }
$$
c_{\lambda}(q,t) \ses 
\prod_{s \in \lambda} \left( 1-q^{a(s)}t^{l(s)+1}\right)
\ssp. \eqno 3.11
$$
For example, up to a scalar,
$$
J_{2,2,1}(X;q,t) \doteq
\left|
\matrix{
S_{2,2,1}[X^{tq}]& S_{3,1, 1}[X^{tq}]& S_{3, 2}[X^{tq}]& S_{4, 1}[X^{tq}]& S_{5}[X^{tq}]\cr\cr
-[1,3,1]& [3,1,1]& 0& 0 & 0\cr\cr 
-[2,0,3]& -[3,0,2]& [3,2, 0]&0& 0\cr\cr 
[1,0,4]& -[0,4, 1]&-[1, 4, 0]&[ 4, 1,0]& 0\cr\cr 
0 & [0, 0, 5]& 0& -[0, 5, 0]& [5,0,0]
}
\right|
\matrix{
{}\cr\cr
{}_{(3,1,1)}\cr\cr
{}_{(3,2,0)}\cr\cr
{}_{(4,1,0)}\cr\cr
{}_{(5,0,0)}
} \eqno 3.12
$$
where an entry $\pm [\alpha]$ stands for 
$\pm \bigl( [|2,2,1|] -[|\alpha|]  \bigr)$.  
Each column of such determinant is obtained from the ordered
expansion of a Schur function.  
For example, the second column is given by the expansion of 
$S_{3,1,1}$ computed in 3.4.
Notice we can limit ourselves to vectors of length 
$\ell(\lambda)$ since each entry is the difference of two vectors of which the last
$n-\ell(\lambda)$ components are equal.
\sa
\noindent {\bf Proof of Theorem 3.1} \quad 
It suffices to check that
$\Bigl([|\lambda|] -M_1 \Bigr) J_{\lambda}'(X;q,t)=0$, where $J_{\lambda}'(X;q,t)$ 
is the determinantal expression 3.9, 
given that the eigenvalues are distinct for each $\lambda$. 
The operator acts only on the first row of the determinant allowing
us to use formula 3.6 for the action of $M_1$ on $S_{\mu}[X^{tq}]$.  
This gives
$$ 
 \bigl( [|\lambda|]-M_1  \bigr)  J'_{\lambda}(X;q,t) \ses
\det \left| \, \matrix{ \sum_{\beta \geq \lambda}  
[\lambda]_{\lambda\beta}  \,
 S^{\beta} & \ldots &  \sum_{\beta \geq \mu} [\lambda]_{\mu\beta} \,
 S^{\beta} & \ldots \cr
                                                \vdots& & \vdots &  \cr
[\lambda]_{\lambda\nu} & \ldots &  
[\lambda]_{\mu\nu} & \ldots \cr
 \vdots & & \vdots &  \cr} \right|\ssp. 
\eqno 3.13
$$
The determinant can be shown to vanish by examining the coefficient
of all $S^\beta$'s.
The coefficient of $S^{\lambda}$  is proportional 
to $[\lambda]_{\lambda\lambda}=0$.  
For $\beta> \lambda$, the coefficient of $S^{\beta}$ is equal to the determinant
 $$  
\det \left| \, \matrix{ [\lambda]_{\lambda\beta} 
  & \ldots &   [\lambda]_{\mu\beta} 
  & \ldots \cr
                                                \vdots& & \vdots &  \cr
 [\lambda]_{\lambda\beta}  & \ldots &  
[\lambda]_{\mu\beta}  & \ldots \cr
\vdots & & \vdots &  \cr} \right|
\matrix{
{\rm row}\,1&\cr
&\cr
&\cr
{\rm row}\,\beta&\cr
&\cr
}
\ssp, \eqno 3.14
$$
which has two identical rows, and
thus also vanishes.
\sas 
To verify the normalization, we recall [M1]
$$
J_{\lambda}(X;q,t) \ses c_{\lambda'}(t,q) \,S_{\lambda}[X^{tq}] 
+ \sum_{\mu > \lambda} c_{\lambda \mu} S_{\mu}[X^{tq}]
\ssp \eqno 3.15
$$
and observe that because, from 3.7,
$[\lambda]_{\mu \nu}=0$ for $\mu>\nu$,
the sub-determinant giving the coefficient of 
$S_{\lambda}[X^{tq}]$ in $J'_\lambda(X;q,t)$ is 
triangular and equals $v_\lambda(q,t)$.\hfill $\square$
\sas
The formula for $J_\lambda(X;q,t)$ in terms of 
$S_\lambda[X^{tq}]$ can be converted into a
similar expression in terms of the $S_{\mu}[X]$ basis
by applying to expression 3.9,
the involution $\omega_{q,t}$ satisfying properties [M1]:
$$
\omega_{q,t} S_{\lambda}[X^{tq}] = S_{\lambda'}[X], 
\qquad \qquad \omega_{q,t} J_{\lambda}(X;q,t) 
\doteq J_{\lambda'}(X;t,q) \,\,\hbox{ up to a scalar}\ssp .
$$
Renormalizing using 
$$
J_{\lambda}(X;q,t)\ses 
c_{\lambda}(q,t)\, S_\lambda[X]+
\sum_{\mu<\lambda}c_{\lambda\mu}\,S_\mu[X]
\ssp, \eqno 3.16
$$
we have
\sa
\noindent {\bf Corollary 3.2}
$$ 
J_{\lambda}(X;q,t) \ses {c_{\lambda}(q,t) \over u_{\lambda}(t,q)} \det \left| \, 
\matrix{ S_{\lambda} [X] & \ldots & S_{\mu}[X] & \ldots \cr
                                                \vdots& & \vdots &  \cr
[\lambda']_{\lambda'\nu'}(t,q) & \ldots &  
 [\lambda']_{\mu'\nu'}(t,q) & \ldots \cr
 \vdots & & \vdots &  \cr} \right|_{\nu< \lambda \atop\mu\leq\lambda}
\ssp,  \eqno 3.17
$$
{\it where $u_{\lambda}(q,t) = \prod_{\mu < \lambda} 
\Bigl( [|\lambda'|]  - [|\mu'|] \Bigr)$. }
\sa
\noindent  The first row is the list of Schur fonctions $S_{\mu}$, $ \mu \leq \lambda$,
the other entries are the polynomials $[\lambda']_{\mu' \nu'}(t,q)$ for $\nu < \lambda$. 
\sap
A determinantal formula for the Macdonald polynomials in terms of
the monomial basis may now be obtained from formula 3.17 .
The interpretation of the entries of the matrix in 3.17 (aside
from the top row) as scalar products will allow that we pass to other 
bases of symmetric functions.  
\sas
We use the space of polynomials in $x_1,\ldots,x_n$ as a 
free module over the ring of symmetric polynomials.  This space
has the scalar product, $ \langle f,g\rangle = \partial_\omega (f g)$\ssp,
and two adjoint bases 
$$
\bigl\{\bar x^{I} = (-x_1)^{I_1} \cdots (-x_m)^{I_m}
\bigr\}_{0^m\subseteq I \subseteq \rho}
 \quad \hbox{and}
\quad 
\bigl\{E_I = e_{I_1}[0] e_{I_{2}}[X_1] \cdots e_{I_m}[X_{m-1}]
\bigr\}_{0^m\subseteq I \subseteq \rho }
\eqno 3.18
$$ 
for $\rho=(0,\dots,m-1)$ and $0^m=(0,\ldots,0)\in \N^m$. That is,
$$
\langle \bar x^I, E_{\rho-J} \rangle = \d_{\omega} \Bigl( \bar x^{I} E_{\rho-J}\Bigr)
= \delta_{IJ}, \qquad \qquad I,J \subseteq \rho\ssp. \eqno 3.19
$$      
This result can be used to determine the coefficients in the following matrix;
$$
S_{\lambda'}[X_n] = {\hbox {det} }\left|
\matrix{ e_{\lambda_1} [X_n]  & e_{\lambda_1+1}[X_{n+1}]  & \cdots & e_{\lambda_1+n-1}
[X_{2n-1}]  \cr
  e_{\lambda_2 -1}[X_n]  & e_{\lambda_2}[X_{n+1}]  & \cdots & e_{\lambda_2 +n-2}[X_{2n-1}]  \cr
  \vdots & \vdots & \ddots & \vdots \cr
  e_{\lambda_{n} -n+1}[X_n]  & e_{\lambda_{n} -n+2}[X_{n+1}]  & \cdots & 
  e_{\lambda_{n}}[X_{2n-1}] \cr} 
  \right| \ssp, \eqno 3.20
$$
which is obtained by applying the $\omega$ involution,
$ \omega S_k[X] =  e_k[X]$, $\omega S_{\lambda}[X]= S_{\lambda'}[X]$,
to 1.3 and increasing the alphabets in each column with 
the relation, $e_k[X_l+x_{l+1}+\cdots+x_{2n-1}] = e_k[X_l] +
\sum_{j=1}^{2n-l-1} e_{k-j}[X_l]\, e_j[x_{l+1}+\cdots+x_{2n-1}]$.
Note that this relation does not change the 
value of the determinant since it corresponds 
to multiplication by a unitriangular matrix.
Now using 3.19 with $m=2n$, the coefficients in the expansion of 
this matrix are
$$S_{\lambda'}[X_n] = \sum_{\alpha' \in \rho'}
 \langle \bar x^{\rho-\alpha'}, S_{\lambda'}[X_n] \rangle \, E_{\alpha'} \ssp, 
\eqno 3.21
 $$ 
where $\alpha'= (0,\dots,0,\alpha_1,\dots,\alpha_n)$ 
and $\rho'=(0,\dots,0,n,\dots,2n-1)$ for $\alpha',\,\rho'\in\N^{2n}$ .
Since the ordered expansion of 3.20, and thus 3.21, must coincide 
with the ordered expansion 3.3 of $S_{\lambda}[X]$ in terms
of complete functions, we have 
$$
\epsilon(\lambda,\alpha) = \langle \bar x^{\rho-\alpha'}, 
S_{\lambda'}[X_n] \rangle\ssp. 
\eqno 3.22
$$
This implies that 3.17 can be rewritten, using 
$$
\langle f \rangle_{\lambda \nu} \ses 
\sum_{\alpha} \langle \bar x^{\rho -\alpha'}, f \rangle 
\, \bigl([|\lambda|]_{t,q}-[|\alpha|]_{t,q} \bigr)\ssp, 
\eqno 3.23
$$
where the sum is over all distinct permutations $\alpha$ of $\nu \in \N^n$, 
as
$$ 
J_{\lambda}(X_n;q,t) \ses {c_{\lambda}(q,t) \over u_{\lambda}(t,q)} \det \left| \, 
\matrix{ S_{\lambda} [X_n] & \ldots & S_{\mu}[X_n] & \ldots \cr
                                                \vdots& & \vdots &  \cr
 \langle S_{\lambda}[X_n] \rangle_{\lambda' \nu'} & \ldots &  
 \langle S_{\mu}[X_n] \rangle_{\lambda' \nu'}  & \ldots \cr
 \vdots & & \vdots &  \cr} \right|_{\mu\leq\lambda\atop\nu<\lambda}
\ssp.  \eqno 3.24
$$
With $\epsilon(\mu,\alpha)$ now defined as a scalar product, we can write Macdonald
polynomials in any linear basis of the space 
generated by $\{ S_{\mu} \}_{\mu \leq \lambda}$.
In particular, we have the following theorem: 
\sa
\noindent {\bf Theorem 3.3} \quad { \it The Macdonald polynomials in terms of monomial symmetric functions are
$$ 
J_{\lambda}(X_n;q,t) \ses {c_{\lambda}(q,t) \over u_{\lambda}(t,q)} \det \left| \, 
\matrix{ m_{\lambda} [X_n] & \ldots & m_{\mu}[X_n] & \ldots \cr
                                                \vdots& & \vdots &  \cr
  \langle m_{\lambda}[X_n] \rangle_{\lambda' \nu'} & \ldots &  
  \langle m_{\mu}[X_n] \rangle_{\lambda' \nu'} & \ldots \cr
 \vdots & & \vdots &  \cr} \right|_{\mu\leq\lambda\atop\nu<\lambda}
\ssp.  \eqno 3.25
$$ }
For example, $J_{2,2}(X;q,t)$ is given by
$$
J_{2,2}(X_4;q,t) = {c_{2,2}(q,t) \over u_{2,2}(t,q)} \det \left| \,
\matrix{ m_{2,2}[X_4] & m_{2,1,1}[X_4] & m_{1,1,1,1}[X_4] \cr
 -[3,1]-[1,3] & [3,1] & 0 \cr
 [4,0]+[0,4] & -3[4,0]-[0,4] & [4,0] \cr}\right| , \eqno 3.26
$$
with $[\alpha]=\bigl( [|2,2|]_{t,q}-[|\alpha|]_{t,q} \bigr)$.  
In this expression, the term $-3[4,0]$ appears, for instance, because
$\langle \bar x^{(0,1,2,3,0,5)},m_{2,1,1}[X_4] \rangle = -3$.  
Again, we work with vectors of length $\ell(\lambda)$, 
since adding zeros to $\alpha$ does not change the result.
We skip the problem of computing efficiently all the scalar products in the matrix.

\vfill\supereject

\noindent{\bol 4. Creation Operators}
\sa
A Macdonald polynomial associated to any partition
can be constructed by repeated application of 
creation operators, $B_k^{(n)}$.  Specifically,
$$
B_k^{(n)}\, J_{\lambda}(X;q,t)
\ses
J_{\lambda+1^k}(X;q,t)\ssp. \eqno 4.1
$$
The creation operators were defined originally 
to be [LV],[KN]
$$
B_k^{(n)}\ses
\sum_{|I|=k}
\sum_{l=0}^k
(-t)^l
x^I \, {
R( X_I,X_I^c/t)\,
\over
R(X_I,X_I^c)
}\,M_l^{(I)}
\ssp, \eqno 4.2
$$
where $I$ is a $k$ subset of $\{1,\ldots,n\}$ and
$M^{(I)}$ acts only on $x_i$ for all $i\in I$.
We will see that the action of these operators on the modified 
Schur function basis may be expressed in a manner similar to 
that of the Macdonald operators on this basis, though
the former increases degrees.  We begin by giving a new expression
for the creation operators. 
\sa
\noindent {\bf Proposition 4.1}\quad
{\it The creation operators acting on the space of symmetric
functions are }
$$ B_k^{(n)}\ses {1\over k_t!}\,
\chi_{(k| n-k)}\,
\chi_{\omega(k)}\,
x_1\cdots x_k\,
(1-t\Omega_k) \cdots (1-t^k\Omega_k)\, 
\eqno 4.3
$$
{\it where } $ \Omega_k= \sigma_{1}\cdots \sigma_{k-1} \Theta_k$.
\sa
\noindent {\bf Proof}\quad
The binomial expansion of the right hand side of expression 4.3 becomes
$$
rhs\ses
{1\over k_t!}
\sum_{l=0}^k
{(-1)^lt^{l+1\choose 2}(t)_k \over (t)_{k-l}(t)_l}\,
\chi_{(k| n-k)}\,
x_1\cdots x_k\, 
\chi_{\omega(k)}\,
\Theta_l\cdots\Theta_1 
\ssp, \eqno 4.4
$$
since $\Omega_k^l=\Theta_l\cdots\Theta_1$ on symmetric functions. 
Further, the image of such a function under $\Omega_k^l$
belongs to $\Sym(l|k-l)$
allowing that we use property 1.12 of $\chi_{\omega(k)}$ to obtain
$$
rhs\ses
{1\over k_t!}
\sum_{l=0}^k
{(-1)^lt^{{l+1\choose 2}+l(k-l) }(t)_k \over (t)_{k-l}(t)_l}\,
l_t!(k-l)_t!
\chi_{(k| n-k)}
x_1\cdots x_k\, 
\chi_{(l|k-l)}\,
\Theta_l\cdots\Theta_1 
\ssp. \eqno 4.5
$$
The definition of the Macdonald operators and
$\chi_{(k|n-k)}$ then imply that
$$
rhs \ses
\sum_{l=0}^k
(-t)^l
\partial_{(k|n-k)}\,
R( X_k,X_k^c/t)\, 
x_1\cdots x_k\, 
M_l^{(k)}
\ssp. \eqno 4.6
$$
Finally, the action of $\partial_{(k|n-k)}$ as described in 1.8
confirms that this expression, and therefore 4.3, 
are in fact equivalent to the original definition for the creation operators. \hfill $\square$
\sa
\noindent{\bf Theorem 4.2} \quad 
{\it The action of $B_k^{(n)}$ on $S_\lambda[X^{tq}]$,
for $\ell(\lambda) \leq k$, can be expressed as 
$$
B_k^{(n)} \, S_{\lambda}[X^{tq}]
\ses
\prod_{i=1}^k
\left(
\left[
\matrix{
\cdots&{0}&{1}&{0}&\cdots\cr
\cdots&{0}&{0}&{0}&\cdots
}
\right]
\,-\,t^{k-i}
\left[
\matrix{
\cdots&{0}&{1}&{0}&\cdots\cr
\cdots&{0}&{X^t}&{0}&\cdots
}
\right]
\right)\, S_{\lambda}[X^{tq}],
\eqno 4.7
$$
where the $i^{th}$ column is the only non-zero column of the operators. }

\noindent 
The relation, $ S_j[X^{tq}+X^t]=S_j[qX^{tq}]=q^jS_j[X^{tq}]\ssp,$
yields immediately
\sa
\noindent{\bf Corollary 4.3}\quad {\it 
Let $f(i,j)=\lambda_{i}+(j-i+1)$.
For any $\ell(\lambda) \leq k$ we have }
$$ 
B_k^{(n)}\,
S_{\lambda}[X^{tq}]
\ses
\left|(1-q^{f(i,j)}t^{k-j})
S_{f(i,j)}[X^{tq}]\right|_{i,j \leq k}
\ssp. \eqno 4.8
$$
For example, when $n=2$ we have
$$
\eqalign{
B_2\,S_{\lambda_1,\lambda_2}[X^{tq}] & \ses
\left(
\left[\matrix{1&0\cr 0&0}\right]-
t \left[\matrix{1&0\cr X^t&0}\right]
\right)
\left(
\left[\matrix{0&1\cr 0&0}\right]-
t^0 \left[\matrix{0&1\cr 0&X^t}\right]
\right)\,S_{\lambda_1,\lambda_2}[X^{tq}]\cr
\cr
&\ses
\left|
\matrix{ (1-q^{\lambda_1+1}t) S_{\lambda_1+1}[X^{tq}] &
(1-q^{\lambda_1+2}) S_{\lambda_1+2}[X^{tq}] \cr
(1-q^{\lambda_2}t) S_{\lambda_2}[X^{tq}] &
(1-q^{\lambda_2+1}) S_{\lambda_2+1}[X^{tq}] 
}
\right|
} $$ 
\sa
\noindent {\bf Proof of Theorem 4.2}\quad 
In what follows, as we have done in the case of the Macdonald operators, 
we shall use $\Theta^t$ instead of $\Theta$ in the operator $B_k$.  This will
allow us to use the action of $M_k$ on $S_{\lambda}[Y+X]$ (with $Y$ an alphabet
invariant under $\Theta^t$) obtained in the previous section.  Letting $Y=X^{tq}-X$
in the final result will prove the theorem.
\sa
\noindent First, notice from 4.3 that $B_k^{(n)}=\chi_{(k|n-k)}B_k^{(k)}$.
This given, we will first describe the action of 
$B_k^{(k)}$ on $S_{\lambda}[Y+X]$.
\sa
\noindent {\bf Claim 4.4} \quad
{\it For $\ell(\lambda) \leq k$ and $\alpha_i\rightarrow\lambda_i-i+1$, 
we have } 
$$ 
\eqalign{
B_k^{(k)}\, 
S_{\lambda}[Y+X]
\ses &
\prod_{i=1}^k
\left(
\left[
\matrix{
\cdots&{0}&{1}&{0}&\cdots\cr
\cdots&{0}&{0}&{0}&\cdots
}
\right]
\,-\,t^{k-i}
\left[
\matrix{
\cdots&{0}&{1}&{0}&\cdots\cr
\cdots&{0}&{X_k^t}&{0}&\cdots
}
\right]
\right) S_{\lambda}[Y+X]
\cr
\ses &
\Big|
G_{\alpha_i+1}(t^{k-1},X_{k})
G_{\alpha_i+2}(t^{k-2},X_{k})
\cdots
G_{\alpha_i+k}(t^{0},X_{k})
\Big|_{i=1..k}
\ssp.
}
\eqno 4.9
$$
{\it where $G_h(t^r,X_k)=S_h[Y+X]-t^r S_h[Y+X+X_k^t]$.  }
\sa
\noindent {\bf Proof of Claim}\quad
With $Y'=Y+X_k^c$, expression 4.2, where $n=k$, yields  
$$ 
B_k^{(k)}\, S_{\lambda}[Y+X]
 = {x_1 \cdots x_k} \sum_{l=0}^k (-t)^l  M_l^{(k)}  
S_{\lambda}[Y'+X_k] , \eqno 4.10
$$
which, using the Macdonald operator action on $S_\lambda[Y'+X_k]$, gives
$$ \eqalign{ B_k^{(k)}\, S_{\lambda}[Y+X] 
& =  {x_1 \cdots x_k}  \sum_{l=0}^k (-t)^l \sum_{1 \leq i_1 < \cdots < i_l \leq k} 
t^{(n-i_1)+ \cdots+ (n-i_l) } \cr
& \qquad \qquad \qquad \times
\left[ \matrix{ 0 & \cdots & 0 &\cdots&0 &\cdots &  0\cr
   0 &\cdots &X_k^t &\cdots&X_k^t& \cdots &0 }\right] 
S_{\lambda}[Y'+X_k] 
\ssp. } \eqno 4.11
$$
In the last expression, for a fixed $l$, $X^t$ occurs in positions $i_1,...,i_l$.
Comparing such terms, we obtain
$$\eqalign{
 & B_k^{(k)}\, S_{\lambda}[Y+X] \cr
& \qquad \ses {x_1 \cdots x_k}  
\prod_{i=1}^k \left( \left[
\matrix{
0&\cdots& 0\cr
0& \cdots & 0 \cr
}
\right]
\,-\,t^{k+1-i}
\left[
\matrix{
0& \cdots & 0 & 0 & 0& \cdots & 0 \cr
0 & \cdots & 0& X_k^t & 0 & \cdots & 0 \cr
}
\right] \right) 
S_{\lambda}[Y+X] 
\ssp ,} \eqno 4.12
$$
or more compactly
$$ B_k^{(k)}\,
S_{\lambda}[Y+X] 
\ses 
x_1 \cdots x_k\,
\Big|
G_{\alpha_i}(t^{k},X_{k})
G_{\alpha_i+1}(t^{k-1},X_{k})
\cdots
G_{\alpha_i+k-1}(t^{1},X_{k})
\Big|_{i=1..k} 
\ssp. \eqno 4.13
$$
The proof of the claim is completed using the following lemma(proved in the appendix):
\sa
\noindent {\bf Lemma 4.5}\quad
{\it For $F_j(t^r,X_k) = S_j[C+X_k]-t^rS_j[C+tX_k]$,
we have
$$ \eqalign{ & \Big|
F_{j+1}  (t^{k-1},X_{k}) F_{j+2}(t^{k-2},X_{k}) 
\ldots
F_{j+k}(1,X_{k})
\Big| 
\qquad \qquad \qquad \qquad \cr
& \qquad \qquad \qquad \ses
x_1\cdots x_k
\Big|
F_{j}(t^{k},X_{k}) F_{j+1}(t^{k-1},X_{k})
\ldots
F_{j+k-1}(t,X_{k})
\Big| \, ,
}
\eqno 4.14 
$$
where the expressions between brackets must be understood as $k \times k$ determinants.}

\noindent In effect, letting
$C \to Y'$, we have $F_j(t^k,X_k) \to G_j(t^k,X_k)$, and thus using 4.14 in 4.13 proves the claim.  \hfill $\square$
\sa
Expression 4.9 gives that the action of 
$B_k^{(n)}=\chi_{(k|n-k)}\,B_k^{(k)}$ on $S_\lambda[Y+X]$ is
$$
B_k^{(n)}\,S_\lambda[Y+X] \ses
\chi_{(k|n-k)}
\Big|G_{\alpha_i+1}(t^{k-1},X_k) \dots G_{\alpha_i+k}(1,X_k)\Big|_{i=1..k}
\ssp.
\eqno 4.15
$$
The proof of Theorem 4.2 is now equivalent to 
showing that $\chi_{(k|n-k)}$ acts on this 
determinant by extending the alphabets from 
$X_k$ to $X$.  This will be achieved with the following 
lemma (proved in
the appendix);
\sa
\noindent {\bf Lemma 4.6}\quad
{\it For any $l\geq k$, we have
$$ \eqalign{ & \Big|G_{j+1}(t^{k-1},X_l) G_{j+2}(t^{k-2},X_l) \dots G_{j+k}(1,X_l)
\Big| \cr
& \qquad = 
\chi_{(1|l-1)}^{(1)}
\cdots
\chi_{(1|l-k+1)}^{(k-1) }
\chi_{(1|l-k)}^{(k)} 
\left|G_{j+1}(t^{k-1},X_{k}) 
G_{j+2}(t^{k-2},X_{k-1}) \dots
G_{j+k}(1,X_{1})\right|
}\ssp,
\eqno 4.16
$$
where, again, the expression between brackets are $k \times k$ determinants. }
\sa

\sa
Application of this lemma to expression 4.15, with $l=k$, gives
$$
\eqalign{ B_k^{(n)}\,& S_\lambda[Y+X] \ses \cr
&
\chi_{(k|n-k)}\,
\chi_{(1|k-1)}^{(1)}
\cdots
\chi_{(1|1)}^{(k-1) }
\Big|G_{\alpha_i+1}(t^{k-1},X_{k})
G_{\alpha_i+2}(t^{k-2},X_{k-1}) \dots
G_{\alpha_i+k}(1,X_{1})\Big|_{i=1,\ldots,k}
}
\ssp. \eqno 4.17
$$
Property 1.14 allows a refactorization of $\chi$ that 
transforms 4.17 into
$$
\eqalign{ B_k^{(n)}\,& S_\lambda[Y+X] \ses \cr
& \chi_{(1|n-1)}^{(1)}
\cdots
\chi_{(1|n-k)}^{(k)} 
\left|G_{\alpha_i+1}(t^{k-1},X_k) 
G_{\alpha_i+2}(t^{k-2},X_{k-1}) \dots
G_{\alpha_i+k}(1,X_1)
\right|_{i=1,\ldots,k}
}
\ssp. \eqno 4.18
$$
Lemma 4.6 may be applied again, with $l=n$, giving
$$
B_k^{(n)} \,
S_{\lambda}[Y+X] 
\ses
\left|G_{\alpha_i+1}(t^{k-1},X)
G_{\alpha_i+2}(t^{k-2},X) \dots
G_{\alpha_i+k}(1,X)
\right|_{i=1,\ldots,k}
\ssp, \eqno 4.19
$$
which proves Theorem 4.2. \hfill $\square$
\sap
\noindent {\bol Appendix}
\sa
\noindent {\bf Proof of Lemma 1.1} \quad 
Using formula 1.2, we have
$$
S_j[X^{tq}+X_k^t]
\ses
S_j[D+x_k^t]
\ses
\sum_{l=0}^j S_{j-l}[D]\,S_l[x_k^t]\,
\ssp,  \eqno A.1
$$
where $D=X^{tq}+X_{k-1}^t$.  Further properties of Schur functions yield, as $x_k^t=tx_k-x_k$, 
$$
\eqalign{
S_j[X^{tq}+X_k^t]
&\ses
S_{j}[D]\,
+
\sum_{l=1}^j (1-1/t)S_{j-l}[D]\,(tx_k)^l\,
\cr
&\ses
S_{j}[D]\,
+
\sum_{l=0}^{j-1} t( x_k-x_k/t)S_{j-l-1}[D]\,(tx_k)^{l}\,
\ssp. 
} \eqno A.2
$$
The alphabet $D$ is invariant under permutations of $x_k, \dots,x_n$,
and thus the first term is invariant up to a constant under $\chi_{(1|n-k)}^{(k)}$.
Thus, using 1.10, we have
$$
\chi_{(1|n-k)}^{(k)}
S_j[X^{tq}+X_k^t]
\ses
(1+t^{-1}+\cdots+t^{-n+k})S_{j}[D]+
\chi_{(1|n-k)}^{(k)}
\sum_{l=0}^{j-1} t(x_k-x_k/t)S_{j-l-1}[D]\,(tx_k)^{l}\,
\ssp, \eqno A.3
$$
and are left only to show
$$
\chi_{(1|n-k)}^{(k)}
\sum_{l=0}^{j-1} t (x_k-x_k/t)S_{j-l-1}[D]\,(tx_k)^{l}\,
\ses
-S_j[D]+S_j[X^{tq}+X^t]\ssp. \eqno A.4
$$
The definition of 
$\chi_{(1|n-k)}^{(k)}$ and
$(x_k-x_k/t)R( x_k,X_k^c/t)= R( x_k,X_{k-1}^c/t)$
convert the left hand side of this expression into
$$
lhs\ses
\d_{n-1}\cdots\d_{k}\,
R( x_k,X_{k-1}^c/t)\,
\sum_{l=0}^{j-1} t S_{j-l-1}[D](tx_k)^{l}
\ssp. \eqno A.5
$$
The identity $S_n[ x-X/t]= R(x,X/t)$  may be used to 
eliminate the resultant, and we obtain
$$
lhs\ses
\d_{n-1}\cdots\d_{k}
\sum_{l=0}^{j-1} t S_{j-l-1}[D](tx_k)^{l}S_{n-k}[ x_k-X_{k-1}^c/t]
\ssp. \eqno A.6
$$
A further Schur function property allows that the factor $x_k^l$ be
used to increase the index of $S_{n-k}$;
$$
lhs
=
\d_{n-1}\cdots\d_{k}
\sum_{l=0}^{j-1} t^{l+1} S_{j-l-1}[D]S_{n-k+l}[x_k-X_{k-1}^c/t]
\ssp. \eqno A.7
$$
Now we can let the divided differences act (see [L1] for similar computations).  Using the property
$\d_k S_l[x_k-X_{k-1}^c/t] =S_{l-1}[x_k+x_{k+1}-X_{k-1}^c/t]$, 
we arrive at the expression,
$$
\eqalign{lhs
=
\sum_{l=0}^{j-1} t^{ l+1} S_{j-l-1}[D]S_{l+1}[x_k+\cdots+x_n-X_{k-1}^c/t]
& =
\sum_{l=0}^{j-1} S_{j-l-1}[D]S_{l+1}[(X_{k-1}^c)^t] \cr
& = -S_j[D] +S_j [D+(X_{k-1}^c)^t] \ssp. \cr} \eqno A.8
$$
Clearly $D+(X_{k-1}^c)^t=X^{tq}+X^t$, which 
implies that expression A.4 holds. \hfill $\square$

\sap
\noindent {\bf Proof of formula 2.11}  \quad We want to show that
$$ \eqalign{ & S_{\lambda}[Y+X] R(X_k,X_k^c)  \cr
& \quad = \Bigl| S_{\alpha_i-k+n} [Y+X_k] \cdots 
S_{\alpha_i-1+n} [Y+X_k] 
S_{\alpha_i+k} [Y+X_k^c] 
\cdots S_{\alpha_i+n-1}[Y+X_k^c] \Bigr|_{i=1..n}\ssp, \cr} \eqno A.9
$$
where $\alpha_i \to \lambda_i -i + 1$.
Let $f_i,\,i=1,\ldots,n$, denote arbitrary one variable polynomials in the space of 
polynomials in $x_1,\dots,x_n$ taken as a free module
over the ring of symmetric polynomials.
From the usual Newton's interpolation in one variable ([L1], Lemma Ni5), 
given such $f_i$, we have
$$
{\hbox {det }} \bigl| \d_{j-1} \dots \d_{1} f_i(x_1)  \bigr|_{1\leq i,j \leq n} \ses 
 {1 \over \Delta(X) } \, \, {\hbox {det }} \bigl| f_i(x_j) \bigr|_{1 \leq i,j \leq n} , \eqno A.10
$$
where $\Delta(X)=\prod_{i<j} (x_j-x_i)$ is the Vandermonde determinant.
If we let $f_i(x) = S_{\alpha_i+n-1}[Y+x]$ in this expression, using
the identity $\d_k S_j[Y+X_k]  = S_{j-1}[ Y+X_{k+1}]$, we obtain
$$
\eqalign{
 {1 \over \Delta(X) } \, \, {\hbox {det }} \Bigl| S_{\alpha_i+n-1}[Y+x_j] \Bigr |_{1 \leq i,j \leq n} & \ses 
 {\hbox {det }} \Bigl| \d_{j-1} \dots \d_{1} S_{\alpha_i+n-1}[Y+x_1]  \Bigr|_{1\leq i,j \leq n} \cr
& \ses
 {\hbox {det }} \Bigl| S_{\alpha_i+n-j}[Y+X_j] \Bigr|_{1\leq i,j \leq n} 
.  } \eqno A.11
$$
The alphabets can be increased in each column of the determinant 
on the right hand side of this expression by using the relation 
$$
S_{\alpha_i+n-j}[Y+X_j]=S_{\alpha_i+n-j}[Y+X] + \sum_{\ell=1}^{n-j} 
S_{\alpha_i+n-j-\ell}[Y+X] S_{\ell}[-X_j^c]\ssp,  
$$
giving that
$$
{\hbox {det }} \Bigl| S_{\alpha_i+n-j}[Y+X_j] \Bigr|_{1\leq i,j \leq n} = {\hbox {det }}
\Bigl| S_{\alpha_i+n-j}[Y+X] \Bigr|_{1\leq i,j \leq n}\ssp. \eqno A.12
$$
The right hand side may thus be substituted into formula
$A.11$ implying that 
$$ 
S_{\lambda}[Y+X] \Delta(X) = \Bigl| S_{\alpha_i+n-1} [Y+x_n] S_{\alpha_i+n-1} [Y+x_{n-1}]
 \cdots S_{\alpha_i+n-1}[Y+x_1] \Bigr|_{i=1..n}\ssp. 
\eqno A.13
$$
This result applied to arbitrary minors of order $n-k$ in the first 
$n-k$ columns and of order $k$ in the last $k$ columns of the 
right hand side of $A.13$, gives by linearity, 
$$ \eqalign{ & S_{\lambda}[Y+X] \Delta(X) \cr
 & \qquad =\Bigl| S_{\alpha_i+n-1} [Y+x_n] \cdots S_{\alpha_i+n-1} [Y+x_{k+1}]
  S_{\alpha_i+n-1} [Y+x_{k}] \cdots S_{\alpha_i+n-1}[Y+x_1] \Bigr|_{i=1..n} \cr
 & \qquad  =\Bigl| S_{\alpha_i+k}[Y+X_k^c] \cdots S_{\alpha_i+n-1}[Y+X_k^c] S_{\alpha_i-k+n}[Y+X_k] 
 \cdots S_{\alpha_i+n-1}[Y+X_k] \Bigr| \cr
 & \qquad \qquad \qquad \qquad \qquad \qquad \qquad \qquad \times \Delta(X_k^c) \Delta (X_k).\cr}
 \eqno A.14$$
Finally, with  $ \Delta(X)/\bigl( \Delta(X_k) 
\Delta(X_k^c)\bigr) = R(X_k,X_k^c) (-1)^{k(n-k)},$ one sees
that $A.14$ is equivalent to $A.9$, which proves the assertion. \hfill $\square$
\sap
\noindent{\bf Proof of Lemma 4.5 }\quad
As a result of $ \d_iS_{j}[A+b_i]= S_{j-1}[A+b_i+b_{i+1}], $
we have
$$
\d_i \,F_j(t^r,A_i)\ses F_{j-1}(t^{r+1},A_{i+1})
\ssp,
\eqno A.15
$$
allowing that the left hand side of 4.14 be rewritten as
$$
lhs\ses
\Big|
\partial_{k-1}..\partial_1 F_{j+k}(1,X_{1}),
\partial_{k-1}..\partial_{2} F_{j+k}(1,X_{2}),
\ldots,
\partial_{k-1} F_{j+k}(1,X_{k-1}),
F_{j+k}(1,X_{k})
\Big|\ssp. \eqno A.16
$$
We can factor out $\d_{k-1} \cdots \d_1$, since all the columns, except the first one are symmetrical
in $x_1,\dots,x_k$.  Similarly we can factor out successively $(\d_{k-1} \cdots \d_2), \dots, \d_{k-1}$
to get
$$
lhs\ses 
(\partial_{k-1}..\partial_{1})
(\partial_{k-1}..\partial_{2})
\cdots
(\partial_{k-1})
\Big|
F_{j+k}(1,X_{1}),\ldots,
F_{j+k}(1,X_{k-1}),
F_{j+k}(1,X_{k})
\Big|
\ssp. \eqno A.17
$$
With $F_j(t^r,X_k) 
= x_k F_{j-1}(t^{r+1} ,X_k) + F_j (t^r,X_{k-1}), $
we can substitute
$x_kF_{j+k-1}(t,X_k)+ F_{j+k}(1,X_{k-1})$
in the last column of the matrix.
This produces a sum of two determinants, of which the
second vanishes thanks to column $k-1$.  
We thus obtain 
$$
lhs\ses 
(\partial_{k-1}..\partial_{1})
(\partial_{k-1}..\partial_{2})
\cdots
(\partial_{k-1})
\Big|
F_{j+k}(t^{0},X_{1}),\ldots,
F_{j+k}(t^0,X_{k-1}),
x_kF_{j+k-1}(t,X_{k})
\Big|
\ssp, \eqno A.18
$$
and by iteration;
$$
lhs=
(\partial_{k-1}..\partial_{1})
\cdots
(\partial_{k-1})
\Big|
x_1F_{j+k-1}(t,X_{1}),\ldots,
x_{k-1}F_{j+k-1}(t,X_{k-1}),
x_k F_{j+k-1}(t,X_{k})
\Big|
\ssp. \eqno A.19
$$
The symmetric factor $x_1\cdots x_k$ can be factored out of the
determinant and the divided differences may be reintroduced;
$$
lhs\ses
x_1\cdots x_k
\Big|
\partial_{k-1}\cdots\partial_{1}
F_{j+k-1}(t,X_{1}),
\ldots,
\partial_{k-1}
F_{j+k-1}(t,X_{k-1}),
F_{j+k-1}(t,X_{k})
\Big|
\ssp. \eqno A.20
$$
Transformation on  columns of this determinant using
relation A.15 completes the proof.
\hfill $\square$
\sap
\noindent {\bf Proof of Lemma 4.6} 
\quad This identity, as it is a consequence of a sucession
of transformations performed on the columns,
can be proven by induction on $k$ given that
$$
G_{j}(t^{r},X_l)\ses
\chi_{(1|l-k)}^{(k)}\,
G_{j}(t^{r},X_{k})
-c(k)\, G_{j}(t^{r},X_{k-1})
\ssp,
\eqno A.21
$$
which is a direct consequence of lemma 1.1. 
If $k=1$, $G_{j+1}(t^0,X_0)=0$ implies that
$$
G_{j+1} (1,X_l)\ses
\chi_{(1|l-1)}^{(1)} \,
G_{j+1} (1,X_1) 
\ssp, \eqno A.22
$$
which is the case $k=1$ of the lemma.
Assuming that the lemma holds for a $k-1$ determinant, we get 
$$ \eqalign{ 
& \Big|G_{j+1}(t^{k-1},X_l) G_{j+2}(t^{k-2},X_l) \dots G_{j+k}(1,X_l)
\Big| 
\qquad \qquad \qquad \qquad
\qquad \qquad \qquad \qquad
\cr
&\qquad 
\ses
\chi_{(1|l-1)}^{(1)}
\cdots
\chi_{(1|l-k+1)}^{(k-1) }
\Big|G_{j+1}(t^{k-1},X_{l}) 
G_{j+2}(t^{k-2},X_{k-1}) \dots
G_{j+k}(1,X_{1})\Big|
\ssp.}
\eqno A.23
$$
The first column becomes
$ \chi_{(1|l-k)}^{(k)} G_{j+1}(t^{k-1},X_{k}) -c(k)\, G_{j+1}(t^{k-1},X_{k-1})$, by 
relation A.21, producing a sum of two determinants.
As the first of these is exactly the desired result,
the second must be shown to vanish.  Specifically,
it must be that 
$$
\chi_{(1|l-1)}^{(1)}
\cdots
\chi_{(1|l-k+1)}^{(k-1) }
\Big|G_{j+1}(t^{k-1},X_{k-1}) 
G_{j+2}(t^{k-2},X_{k-1}) \dots
G_{j+k}(1,X_{1})\Big| 
=0\ssp. \eqno A.24
$$
The left hand side of this expression may be transformed 
first by using property 1.14 to refactor
$ \chi_{(1|l-1)}^{(1)} \cdots \chi_{(1|l-k+1)}^{(k-1) } $ 
into the form
$ \chi_{(k-1|n-k+1)} \chi_{(1|k-2)}^{(1)} 
\cdots \chi_{(1|1)}^{(k-2) } $,
and then by using the induction hypothesis, where $l=k-1$; 
$$ \eqalign{ 
lhs\ses
& 
 \chi_{(k-1|n-k+1)} \chi_{(1|k-2)}^{(1)}
\cdots
\chi_{(1|1)}^{(k-2) }
\Big|G_{j+1}(t^{k-1},X_{k-1}) 
G_{j+2}(t^{k-2},X_{k-1}) \dots
G_{j+k}(1,X_{1})\Big| 
\cr
\ses&
 \chi_{(k-1|n-k+1)} \Big|G_{j+1}(t^{k-1},X_{k-1}) 
G_{j+2}(t^{k-2},X_{k-1}) \dots G_{j+k}(1,X_{k-1})\Big| 
\ssp.} \eqno A.25
$$
We now have $ \chi_{(k-1|n-k+1)}$ applied to a determinant with $k$ columns that 
satisfies the special case of $x_k=0$ in Lemma 4.5.
Thus, the determinant vanishes, proving the lemma.
\hfill $\square$
\sap
\noindent {\bol Acknowledgments}
\sa
The author, J. Morse, expresses thanks to Adriano Garsia for providing 
the (NSF) support that made possible this joint work. 
L. Lapointe is supported through an NSERC postdoctoral fellowship.
We have extensively used the algebraic combinatorics environment, ACE,
Maple library [V]. 
\sap
\noindent {\bol References}
\sa
\item{[GH]} A.M. Garsia and M. Haiman, {\it A graded representation module for Macdonald's
polynomials}, Proc. Natl. Acad. Sci. USA V {\bf 90} (1993) 3607--3610.
\smallskip
\item{[H]} F. Hirzebruch, {\it Topological methods in algebraic geometry}, 3rd ed., Springer, 1966.
\smallskip
\item{[HBJ]} F. Hirzebruch, T. Berger and R. Jung, {\it Manifolds and Modular Forms}, Aspects of Mathematics,
Vieweg, Braunschweig/Wiesbaden, 1992.
\smallskip
\item{[KN]} A. Kirillov and M. Noumi, { \it Affine Hecke algebras and raising operators for Macdonald
polynomials}, Duke Math. J., to appear.
\smallskip
\item{[L1]} A. Lascoux, {\it Notes on interpolation in one and several variables}, 
http://phalanstere.univ-mlv.fr/{\~{}}al/.
\smallskip
\item{[L2]} A. Lascoux, {\it Fonctions de Schur et grasmanniennes}, C.R. Acad. Sc. Paris {\bf 281}
 S\'erie A (1975),
813--815.
\smallskip
\item{[LS]} A. Lascoux and M.-P. Sch\"utzenberger, {\it Polyn\^omes de Schubert}, 
C.R. Acad. Sc. Paris {\bf 294} S\'erie
I (1982), 447--450.
\smallskip
\item{[LV]} L. Lapointe and L. Vinet, {A short proof of the integrality of the Macdonald $(q,t)$-Kostka coefficients},
Duke Math. J. {\bf 91} (1998), 205--214. 
\smallskip
\item{[M1]} I.G. Macdonald, {\it Symmetric Functions and Hall Polynomials}, 2d ed., Oxford Math. Monographs,
Clarendon Press, New York, 1995.
\smallskip
\item{[M2]} I.G. Macdonald, {\it Notes on Schubert Polynomials}, Publ. du LACIM, Montr\'eal, 1991.
\smallskip
\item{[RS]} S.N.M. Ruijsenaars and H. Schneider, { \it A new class of integrable systems and its relation to 
solitons}, Ann. Phys.  {\bf 170} (1986), 370--405.
\smallskip
\item{[V]} {S. Veigneau}, {\it ACE, an Algebraic
                      Combinatorics Environment for the computer
                      algebra system MAPLE\/}, {\it Version 3.0\/},
                      Universit\'e de Marne-la-Vall\'ee, 1998,
http://phalanstere.
univ-mlv.fr/{\~{}}ace/.

\end